\documentclass[12pt,thmsa]{article}%
\usepackage{sw20bams}
\usepackage{amsmath}
\usepackage{amsfonts}
\usepackage{amssymb}
\usepackage{graphicx}%
\providecommand{\U}[1]{\protect\rule{.1in}{.1in}}
\begin{document}

\author{Steven R. Finch and Li-Yan Zhu}
\title{Searching for a Shoreline}
\date{March 12, 2016}
\maketitle

\begin{abstract}
Logarithmic spirals are conjectured to be optimal escape paths from a half
plane ocean. Assuming this, we find the rate of increase for both min-max and
min-mean interpretations of ``optimal''. For the one-dimensional analog, which
we call logarithmic coils, our min-mean solution differs from a widely-cited
published account.

\end{abstract}

\footnotetext{Copyright \copyright \ 2005, 2016 by Steven R. Finch. All rights
reserved.}A ship is lost in a dense fog at sea and must reach land as soon as
possible. The captain knows that the shore is straight line, but has no
information about its distance or direction. Equivalently, the sea is known to
be a half plane, but the ship's location and orientation relative to the
boundary is unknown. Assuming its speed is constant, what is the best path for
the ship to follow in its search for the shore?

The word \textquotedblleft best\textquotedblright\ can be understood in
several ways \cite{FW}. We start with minimizing the worst-case scenario
(min-max); a relevant conjecture that the family of logarithmic spirals
contains the minimal path remains open. Our small contribution is that of
providing the computational details that underlie a proposition due to
Baeza-Yates, Culberson \&\ Rawlins \cite{BY1, BY2, BY3}. We then adopt a
different sense of \textquotedblleft best\textquotedblright\ and determine the
logarithmic spiral that minimizes the expected pathlength (min-mean), in which
shoreline directions are assumed to be uniformly distributed. Except for the
(admittedly large) theoretical gap regarding the optimality of logarithmic
spirals, the calculations in this two-dimensional setting are straightforward.

We subsequently turn to the one-dimensional analog of the search problem. The
shore is now simply a point on a line and the spirals here are necessarily
self-intersecting. A\ large computer science literature on this problem
exists. The solution of the min-max problem was first found by Beck \&\ Newman
\cite{BN}. Their approach to the min-mean problem, however, suffers from the
assumption of a nonuniform target distribution (a certain scaling property,
true in the two-dimensional setting, is less apparent here). We give our
solution, which is distinct from theirs, and hope to initiate discussion on
this issue.

The three-dimensional analog, for which shores are planes in space, would seem
to be very difficult. We wonder if an appropriate extension of spiral has ever
been examined in the past.

\subsection{Planar Setting:\ Min-Max}

Let $\kappa>0$. Three preliminary results are:\medskip

\noindent\textbf{Lemma 1.} \textit{The distance between the line} $Ax+By+C=0$
\textit{and the origin is} $|C|/\sqrt{A^{2}+B^{2}}$.\medskip

\noindent\textbf{Lemma 2.} \textit{The equation of a line tangent to the
circle of radius} $R$, \textit{center at the origin, is} $r=R\sec
(\theta-\omega)$\textit{, where }$\omega$ \textit{corresponds to the point of
tangency.}\medskip

\noindent\textbf{Lemma 3.} \textit{The equation of a line tangent to the
spiral} $r=e^{\kappa\theta}$ \textit{is} $y-e^{\kappa\theta}\sin
(\theta)=m(x-e^{\kappa\theta}\cos(\theta))$\textit{, where }$\theta$\textit{
corresponds to the point of tangency and the slope is given by}
\[
m=\frac{\kappa\sin(\theta)+\cos(\theta)}{\kappa\cos(\theta)-\sin(\theta)}.
\]

\noindent\textbf{Proof of Lemma 1.} The unit vector $(A,B)/\sqrt{A^{2}+B^{2}}$
is normal to the line $Ax+By+C=0$, hence the point $(-CA,-CB)/(A^{2}+B^{2})$
on the line determines its distance from $(0,0)$.\medskip

\noindent\textbf{Proof of Lemma 2.} In rectangular coordinates, the line is
given by
\[
y=R\sin(\omega)-\cot(\omega)\left(  x-R\cos(\omega)\right)  .
\]
In polar coordinates, therefore, we have
\[
r\sin(\theta)=R\sin(\omega)-\cot(\omega)\left(  r\cos(\theta)-R\cos
(\omega)\right)
\]
and so
\[
\frac rR=\frac{\sin(\omega)+\cot(\omega)\cos(\omega)}{\sin(\theta)+\cot
(\omega)\cos(\theta)}=\sec(\theta-\omega).
\]

\noindent\textbf{Proof of Lemma 3.} Clearly
\[
\frac{dy}{dx}=\frac{dy/d\theta}{dx/d\theta}=\frac{(e^{\kappa\theta}\sin
(\theta))^{\prime}}{(e^{\kappa\theta}\cos(\theta))^{\prime}}=\frac{\kappa
\sin(\theta)+\cos(\theta)}{\kappa\cos(\theta)-\sin(\theta)}.
\]

\noindent\textbf{Theorem 4.} \textit{Of all lines tangent to the spiral}
$r=e^{\kappa\theta}$\textit{, there is exactly one that is tangent to the
circle of radius} $R$\textit{, center at the origin. Call this line}
$L$\textit{. The tangency point of} $L$ \textit{with the spiral} \textit{is}
\[
\theta_{0}=\frac{1}{\kappa}\left(  \ln(R)+\frac{1}{2}\ln(1+\kappa^{2})\right)
.
\]
\textit{The tangency point of }$L$ \textit{with the circle is}
\[
\omega_{0}=\theta_{0}-\arccos\left(  \frac{1}{\sqrt{1+\kappa^{2}}}\right)
<\theta_{0};
\]
\textit{thus} $L$ \textit{has equation} $r=R\sec(\theta-\omega_{0})$.\medskip

\noindent\textbf{Proof of Theorem 4.} Apply Lemma 1 with $A=m$, $B=-1$ and
$C=e^{\kappa\theta}(\sin(\theta)-m\cos(\theta))$ to obtain $R^{2}$ as an
expression in $\kappa$, $\theta$, $m$. Lemma 3 further gives $m$ as an
expression in $\kappa$, $\theta$. We find that $R^{2}(1+\kappa^{2}%
)=e^{2\kappa\theta}$ and hence the formula for $\theta_{0}$ is true. By Lemma
2, $e^{\kappa\theta_{0}}=R\sec(\theta_{0}-\omega)$ and thus the formula for
$\omega_{0}$ is true.\newline

Think of a ship, starting at the origin and moving along the spiral. Its first
contact point with $L$ is at $\theta_{0}$. We wish to compute its second
contact point $\theta_{1}$. The reason is that, in the interval $\theta
_{0}<\theta<\theta_{1}$, the spiral intersects all other tangent lines to the
circle of radius $R$. At $\theta=\theta_{1}$, repetition begins so we stop
there:\ All possible shorelines at distance $R\,$ from the origin have at this
point been found.

A\ closed-form expression for $\theta_{1}$ is not known, but it uniquely
satisfies the equation
\[%
\begin{array}
[c]{ccc}%
e^{\kappa\theta_{1}}=R\sec(\theta_{1}-\omega_{0}), &  & \theta_{0}<\theta
_{1}<\theta_{0}+2\pi\text{.}%
\end{array}
\]
Once we have $\theta_{1}$ for $R=1$, we have it for all $R$ via the formula
\[
\theta_{1}(R)=\frac{1}{\kappa}\ln(R)+\theta_{1}(1)
\]
since
\[
\omega_{0}(R)=\frac{1}{\kappa}\ln(R)+\omega_{0}(1),
\]
so $\theta_{1}(R)-\omega_{0}(R)=\theta_{1}(1)-\omega_{0}(1)$ and thus
\[
e^{\kappa\theta_{1}(R)}=Re^{\kappa\theta_{1}(1)}=R\sec(\theta_{1}%
(1)-\omega_{0}(1))=R\sec(\theta_{1}(R)-\omega_{0}(R)).
\]
Such scaling behavior is valuable here -- we may consider $R=1$ without loss
of generality -- but this property fails in Section 3.\medskip

\noindent\textbf{Lemma 5.} \textit{The arclength of the spiral} $r=e^{\kappa
\theta}$ \textit{up to} $\Theta$ \textit{is}
\[
\sqrt{1+\kappa^{2}}%
%TCIMACRO{\dint \limits_{-\infty}^{\Theta}}%
%BeginExpansion
{\displaystyle\int\limits_{-\infty}^{\Theta}}
%EndExpansion
e^{\kappa\theta}d\theta=\frac{\sqrt{1+\kappa^{2}}}{\kappa}e^{\kappa\Theta
}\text{.}%
\]

\noindent\textbf{Proof of Lemma 5.} From $dr=\kappa\,e^{\kappa\theta}d\theta$,
we deduce that $ds^{2}=r^{2}d\theta^{2}+dr^{2}=e^{2\kappa\theta}d\theta
^{2}+\kappa^{2}e^{2\kappa\theta}d\theta^{2}=(1+\kappa^{2})e^{2\kappa\theta
}d\theta^{2}$.\medskip

With the assumption that $R=1$, the two-dimensional min-max problem reduces to
minimizing $\left(  \sqrt{1+\kappa^{2}}/\kappa\right)  e^{\kappa\theta_{1}}$
as a function of $\kappa$. While an explicit formula for $\theta_{1}$ in terms
of $\kappa$ is unavailable, a purely numerical scheme suffices to give
$\kappa=0.2124695594...$ with arclength $13.8111351795...$. The latter is
consistent with the estimate $13.81$ reported in \cite{BY1}; earlier estimates
$0.22325$ and $13.49$ from \cite{BY2, BY3} arose when erroneously minimizing
$e^{\kappa\theta_{1}}/\kappa$.

We now obtain trigonometric equations that serve to define the best spiral
more precisely.%
%TCIMACRO{\FRAME{ftbpFU}{5.3255in}{4.2722in}{0pt}{\Qcb{Two helpful pictures for
%the proof of Theorems 6 and 7.}}{}{spiral.eps}%
%width 5.3255in;  height 4.2722in;  depth 0pt;  original-width 7.6544in;
%original-height 6.1324in;  cropleft "0";  croptop "1";  cropright "1";
%cropbottom "0";  filename 'spiral.eps';file-properties "XNPEU";}}}%
%BeginExpansion
\begin{figure}[ptb]%
\centering
\includegraphics[
height=4.2722in,
width=5.3255in
]%
{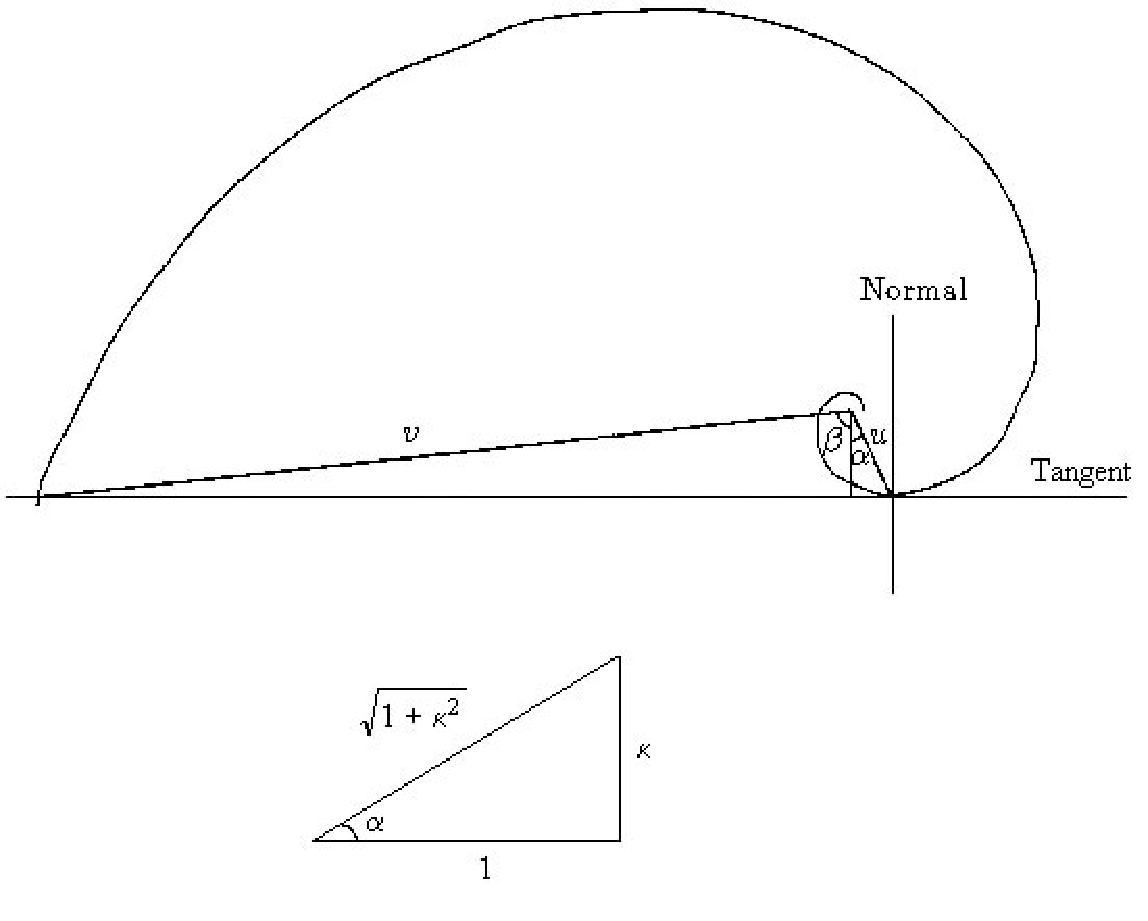}%
\caption{Two helpful pictures for the proof of Theorems 6 and 7.}%
\end{figure}
%EndExpansion
\medskip

\noindent\textbf{Theorem 6.} \textit{The min-max logarithmic spiral has
parameter} $\kappa=\tan\alpha=0.2124695594...=\ln(1.2367284662...)$
\textit{with arclength} $\csc\alpha\sec\beta=13.8111351795...$\textit{, where}
$\alpha$, $\beta$ \textit{satisfy the simultaneous equations}
\[%
\begin{array}
[c]{ccc}%
\dfrac{1}{\tan\alpha}+\dfrac{1}{\tan\beta}=\dfrac{2\pi-\alpha-\beta}{\cos
^{2}\alpha}, &  & \dfrac{\cos\alpha}{\cos\beta}=e^{(2\pi-\alpha-\beta
)\tan\alpha}.
\end{array}
\]

\noindent\textbf{Proof of Theorem 6.} Define angles $\alpha$, $\beta$ and
lengths $u$, $v$ by
\[%
\begin{array}
[c]{lll}%
\theta_{0}=\alpha+\omega_{0}, &  & u=e^{\kappa\theta_{0}}=\sec\alpha,\\
\theta_{1}=(2\pi-\alpha-\beta)+\theta_{0}, &  & v=e^{\kappa\theta_{1}}%
=\sec\beta.
\end{array}
\]
Differentiating with respect to $\alpha$, we obtain
\begin{equation}
u^{\prime}=\sec\alpha\tan\alpha=u\tan\alpha, \label{up}%
\end{equation}
\[
v^{\prime}=\beta^{\prime}\sec\beta\tan\beta=\beta^{\prime}v\tan\beta,
\]
that is,
\begin{equation}
\beta^{\prime}=\frac{v^{\prime}}v\cot\beta=v^{\prime}\cos\beta\cot\beta.
\label{bp}%
\end{equation}
From $e^{\kappa\theta_{1}}=e^{\kappa\theta_{0}}e^{\kappa(2\pi-\alpha-\beta)}$,
it follows that
\begin{equation}
v=u\,e^{(2\pi-\alpha-\beta)\tan\alpha} \label{vu}%
\end{equation}
hence
\begin{align}
v^{\prime}  &  =\left(  (2\pi-\alpha-\beta)\sec^{2}\alpha-(1+\beta^{\prime
})\tan\alpha\right)  v+u^{\prime}\frac vu\label{vp1}\\
&  =\left(  -v^{\prime}\tan\alpha\cos\beta\cot\beta+(2\pi-\alpha-\beta
)\sec^{2}\alpha\right)  v\nonumber
\end{align}
by \label{up}(1) and \label{bp}(2). Since the objective function
\[
e^{\kappa\theta_{1}}\frac{\sqrt{1+\kappa^{2}}}\kappa=v\csc\alpha
\]
is minimized when
\begin{equation}
v^{\prime}\csc\alpha-v\csc\alpha\cot\alpha=0, \label{vp2}%
\end{equation}
we have the additional formula $v^{\prime}=v\cot\alpha$. Substituting this
twice into \label{vp1}(4) yields
\[
\cot\alpha=-v\cos\beta\cot\beta+(2\pi-\alpha-\beta)\sec^{2}\alpha,
\]
thus
\[
\cot\alpha+\cot\beta=(2\pi-\alpha-\beta)\sec^{2}\alpha.
\]
Also, \label{vu}(3) implies immediately that
\[
\sec\beta=\sec\alpha\;e^{(2\pi-\alpha-\beta)\tan\alpha}
\]
which we observe is true independent of \label{vp2}(5).

\subsection{Planar Setting:\ Min-Mean}

On the basis of Section 1, the two-dimensional min-mean problem clearly
reduces to minimizing the average arclength
\begin{align*}
&  \frac{1}{2\pi}\frac{\sqrt{1+\kappa^{2}}}{\kappa}%
%TCIMACRO{\dint \limits_{\omega_{0}}^{\omega_{0}+2\pi}}%
%BeginExpansion
{\displaystyle\int\limits_{\omega_{0}}^{\omega_{0}+2\pi}}
%EndExpansion
e^{\kappa\theta}d\omega\\
&  =\frac{\sqrt{1+\kappa^{2}}}{2\pi\kappa}\left[
%TCIMACRO{\dint \limits_{\theta_{0}}^{0}}%
%BeginExpansion
{\displaystyle\int\limits_{\theta_{0}}^{0}}
%EndExpansion
e^{\kappa\theta}\left(  1-\frac{\kappa}{\sqrt{e^{2\kappa\theta}-1}}\right)
d\theta+%
%TCIMACRO{\dint \limits_{0}^{\theta_{1}}}%
%BeginExpansion
{\displaystyle\int\limits_{0}^{\theta_{1}}}
%EndExpansion
e^{\kappa\theta}\left(  1+\frac{\kappa}{\sqrt{e^{2\kappa\theta}-1}}\right)
d\theta\right] \\
&  =\frac{\sqrt{1+\kappa^{2}}}{2\pi\kappa}\left[  \frac{e^{\kappa\theta_{1}}%
}{\kappa}+\ln\left(  e^{\kappa\theta_{1}}+\sqrt{e^{2\kappa\theta_{1}}%
-1}\right)  -\frac{e^{\kappa\theta_{0}}}{\kappa}+\ln\left(  e^{\kappa
\theta_{0}}+\sqrt{e^{2\kappa\theta_{0}}-1}\right)  \right]
\end{align*}
as a function of $\kappa$, assuming $R=1$. Define for convenience
\[
\Phi(\alpha,\beta)=\left(  -2\csc\alpha+\ln(\sec\alpha+\tan\alpha)+\ln
(\sec\beta+\tan\beta)\right)  \left(  \cot\alpha+\cot\beta\right)  ,
\]%
\[
\Psi(\alpha,\beta)=(\alpha+\beta-2\pi)(\sec\alpha\csc\beta+\csc\alpha\sec
\beta)\sec\alpha,
\]%
\[
\Xi(\alpha,\beta)=\sec\alpha-\cot\alpha\csc\beta+(\tan\alpha\cot\beta
-\csc\alpha\csc\beta)\sec\alpha-(\cot^{2}\alpha+\csc^{2}\alpha)\sec
\beta.\medskip
\]
\noindent\textbf{Theorem 7}. \textit{The min-mean logarithmic spiral has
parameter} $\kappa=\tan\alpha=0.3732051316...=\ln(1.4523822387...)$
\textit{with arclength}
\[
\frac{1}{2\pi}\left(  \ln(\sec\alpha+\tan\alpha)+\ln(\sec\beta+\tan
\beta)-(\sec\alpha-\sec\beta)\cot\alpha\right)  \csc\alpha=7.0321857865...,
\]
\textit{where} $\alpha$, $\beta$ \textit{satisfy the simultaneous equations}
\[%
\begin{array}
[c]{ccc}%
\Phi(\alpha,\beta)+\Psi(\alpha,\beta)=\Xi(\alpha,\beta), &  & \dfrac
{\cos\alpha}{\cos\beta}=e^{(2\pi-\alpha-\beta)\tan\alpha}.
\end{array}
\]

\noindent\textbf{Proof of Theorem 7.} By the observation at the end of the
proof of Theorem 6, the truth of the second equation is independent of the
objective function. It remains to derive the first equation. Define
\[
w=(v-u)\cot\alpha+\ln\left(  v+\sqrt{v^{2}-1}\right)  +\ln\left(
u+\sqrt{u^{2}-1}\right)  ,
\]
then
\begin{align}
w^{\prime}  &  =(v^{\prime}-u^{\prime})\cot\alpha-(v-u)\csc^{2}\alpha
+\frac{v^{\prime}}{\sqrt{v^{2}-1}}+\frac{u^{\prime}}{\sqrt{u^{2}-1}%
}\label{wp1}\\
&  =v^{\prime}\cot\alpha-u-(v-u)\csc^{2}\alpha+\frac{v^{\prime}}{\sqrt
{v^{2}-1}}+\frac{u\tan\alpha}{\sqrt{u^{2}-1}}\nonumber
\end{align}
using \label{up}(1). The objective function
\[
\frac{\sqrt{1+\kappa^{2}}}{2\pi\kappa}\left[  \frac{e^{\kappa\theta_{1}}%
}\kappa+\ln\left(  e^{\kappa\theta_{1}}+\sqrt{e^{2\kappa\theta_{1}}-1}\right)
-\frac{e^{\kappa\theta_{0}}}\kappa+\ln\left(  e^{\kappa\theta_{0}}%
+\sqrt{e^{2\kappa\theta_{0}}-1}\right)  \right]  =\frac{w\csc\alpha}{2\pi}
\]
is minimized when $w^{\prime}=w\cot\alpha$, as with \label{vp2}(5). Between
\label{wp1}(6) and this additional formula, we eliminate $w^{\prime} $ and
solve for $v^{\prime}$ in terms of $\alpha$, $\beta$, $u$, $v$. Substituting
the resulting expression for $v^{\prime}$ into \label{vp1}(4) gives an
equation involving $\Phi$, $\Psi$ and $\Xi$.

\subsection{Linear Setting:\ Min-Max}

Let $\gamma>1$. The one-dimensional analog of the logarithmic spiral
$r=e^{\kappa\theta}$ we study here is
\[
x=(-\gamma)^{\left\lfloor t\right\rfloor }\left(  1-(\gamma+1)(t-\left\lfloor
t\right\rfloor )\right)  .
\]
For lack of standard phraseology, we call this a \textbf{logarithmic coil}.
Local maximum points occur at $(x,t)=(\gamma^{2i},2i)$ where $i$ is an
integer; local minimum points occur at $(x,t)=(-\gamma^{2i-1},2i-1)$.

Given a point $X>0$, the distance $\delta$ that the ship travels to reach $X$
is
\[
\delta=X+2%
%TCIMACRO{\dsum \limits_{j=-\infty}^{2i+1}}%
%BeginExpansion
{\displaystyle\sum\limits_{j=-\infty}^{2i+1}}
%EndExpansion
\gamma^{j}=X+\frac{2\gamma^{2i+2}}{\gamma-1}%
\]
where $\gamma^{2i}<X\leq\gamma^{2i+2}$, that is,
\[
i=\left\lceil \frac{\ln(X)}{2\ln(\gamma)}-1\right\rceil =\left\lfloor \frac
{1}{2}\left\lceil \frac{\ln(X)}{\ln(\gamma)}-1\right\rceil \right\rfloor .
\]
Given a point $X<0$, the corresponding distance $\delta$ is
\[
\delta=-X+2%
%TCIMACRO{\dsum \limits_{j=-\infty}^{2i}}%
%BeginExpansion
{\displaystyle\sum\limits_{j=-\infty}^{2i}}
%EndExpansion
\gamma^{j}=-X+\frac{2\gamma^{2i+1}}{\gamma-1}%
\]
where $\gamma^{2i-1}<-X\leq\gamma^{2i+1}$, that is,
\[
i=\left\lceil \frac{\ln(-X)}{2\ln(\gamma)}-\frac{1}{2}\right\rceil
=\left\lceil \frac{1}{2}\left\lceil \frac{\ln(-X)}{\ln(\gamma)}-1\right\rceil
\right\rceil .
\]
Scaling as observed in Section 1 no longer works here:\ for the points $\pm1$,
we have
\[%
\begin{array}
[c]{ccc}%
\delta(1)=1+\dfrac{2}{\gamma-1}, &  & \delta(-1)=1+\dfrac{2\gamma}{\gamma-1}%
\end{array}
\]
which are not easily related to $\delta(X)$ and $\delta(-X)$. An analysis of
$\delta(X)$ and $\delta(-X)$, which possess sizeable jump discontinuities at
$\gamma^{2i}$ and $-\gamma^{2i-1}$, would seem to require different tools than
before. Computer scientists traditionally normalize by $|X|$; see Figure 2 for
a sample result. For more on the following theorem, see \cite{BN, Fra, Gal,
AG, JS, LOS, HIKL, Sch}.%
%TCIMACRO{\FRAME{ftbpFU}{6.0182in}{2.2174in}{0pt}{\Qcb{Graph of $\delta(X)/|X|$
%for $\gamma=2.$}}{}{mas.eps}{\special{ language "Scientific Word";
%type "GRAPHIC";  maintain-aspect-ratio TRUE;  display "USEDEF";
%valid_file "F";  width 6.0182in;  height 2.2174in;  depth 0pt;
%original-width 7.6268in;  original-height 2.7916in;  cropleft "0";
%filename 'mas.eps';file-properties "XNPEU";}}}%
%BeginExpansion
\begin{figure}[ptb]%
\centering
\includegraphics[
height=2.2174in,
width=6.0182in
]%
{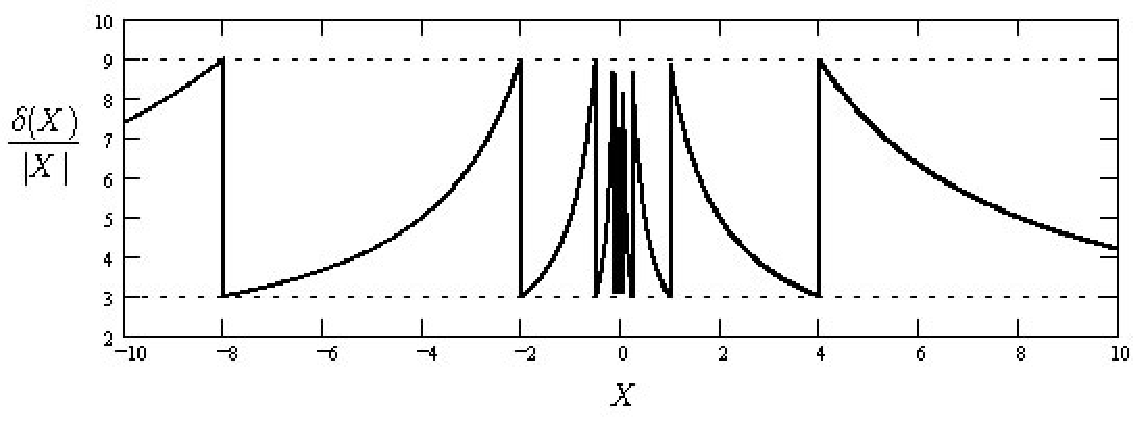}%
\caption{Graph of $\delta(X)/|X|$ for $\gamma=2.$}%
\end{figure}
%EndExpansion
\medskip

\noindent\textbf{Theorem 8. }\textit{The min-max logarithmic coil has
parameter} $\gamma=2$ \textit{with worst-case ratio} $\delta/|X|=9$.\medskip

\noindent\textbf{Proof of Theorem 8.} For simplicity, we examine only positive
$X$. If $X=\gamma^{2k+\varepsilon}$ for some small $\varepsilon>0 $, then
$i=k$ and $\delta/X\rightarrow(2\gamma^{2}+\gamma-1)/(\gamma-1) $ as
$\varepsilon\rightarrow0^{+}$. Calculus gives that $\gamma=2$ is the critical
point, which yields in turn the least maximum value $\delta/|X|=9$.\newline

\subsection{Linear Setting:\ Min-Mean}

The use of $\delta/|X|$ in defining the min-max coil in Section 3 seems fairly
natural; the formulation behind a min-mean coil, however, requires some
careful thought. Consider the integral
\[
I(X)=\frac{1}{2X}%
%TCIMACRO{\dint \limits_{-X}^{X}}%
%BeginExpansion
{\displaystyle\int\limits_{-X}^{X}}
%EndExpansion
\frac{\delta(x)}{|x|}\,dx
\]
whose graph appears in Figure 3. (The minimum and maximum values are
$1+6\ln(2)=5.1588...$ and $1+12/e=5.4145...$ when $\gamma=2$.)
%TCIMACRO{\FRAME{ftbpFU}{6.0857in}{2.4742in}{0pt}{\Qcb{Graph of $I(X)$ for
%$\gamma=2.$}}{}{mas3.eps}{\special{ language "Scientific Word";
%type "GRAPHIC";  maintain-aspect-ratio TRUE;  display "USEDEF";
%valid_file "F";  width 6.0857in;  height 2.4742in;  depth 0pt;
%original-width 7.0828in;  original-height 2.8617in;  cropleft "0";
%filename 'mas3.eps';file-properties "XNPEU";}}}%
%BeginExpansion
\begin{figure}[ptb]%
\centering
\includegraphics[
height=2.4742in,
width=6.0857in
]%
{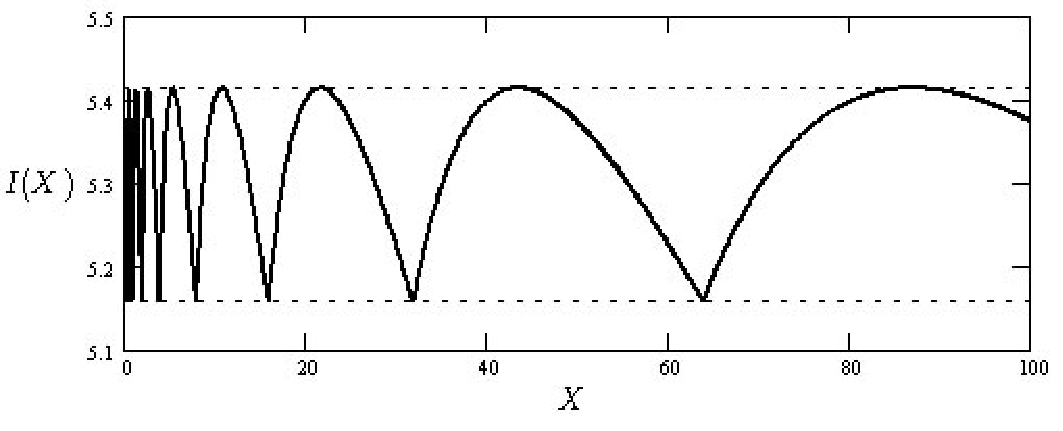}%
\caption{Graph of $I(X)$ for $\gamma=2.$}%
\end{figure}
%EndExpansion
We wish to determine $\gamma$ for which $I$, a kind of normalized average, is
minimal. Since $I$ itself is a periodic function of $X$ (although smoother
than $\delta/|X|$), the word \textquotedblleft minimal\textquotedblright\ can
be used only loosely.

Assuming $\gamma^{2i}<X\leq\gamma^{2i+2}$ and $-\infty<p<i$, we have
\[%
\begin{array}
[c]{ccc}%
%TCIMACRO{\dint \limits_{\gamma^{2p}}^{\gamma^{2p+2}}}%
%BeginExpansion
{\displaystyle\int\limits_{\gamma^{2p}}^{\gamma^{2p+2}}}
%EndExpansion
\dfrac1x\left(  x+\dfrac{2\gamma^{2p+2}}{\gamma-1}\right)  dx=f(p,\gamma
,\gamma^{2p+2}), &  &
%TCIMACRO{\dint \limits_{\gamma^{2i}}^{X}}%
%BeginExpansion
{\displaystyle\int\limits_{\gamma^{2i}}^{X}}
%EndExpansion
\dfrac1x\left(  x+\dfrac{2\gamma^{2i+2}}{\gamma-1}\right)  dx=f(i,\gamma,X)
\end{array}
\]
where
\[
f(i,\gamma,X)=\left(  X-\gamma^{2i}\right)  +\frac{2\gamma^{2i+2}}{\gamma
-1}\left(  \ln(X)-2i\ln(\gamma)\right)  .
\]
Assuming $-\gamma^{2j+1}\leq X<-\gamma^{2j-1}$ and $-\infty<q<j$, we have
\[%
\begin{array}
[c]{ccc}%
%TCIMACRO{\dint \limits_{-\gamma^{2q+1}}^{-\gamma^{2q-1}}}%
%BeginExpansion
{\displaystyle\int\limits_{-\gamma^{2q+1}}^{-\gamma^{2q-1}}}
%EndExpansion
\dfrac1x\left(  -x+\dfrac{2\gamma^{2q+1}}{\gamma-1}\right)  dx=g(q,\gamma
,-\gamma^{2q+1}), &  &
%TCIMACRO{\dint \limits_{X}^{-\gamma^{2j-1}}}%
%BeginExpansion
{\displaystyle\int\limits_{X}^{-\gamma^{2j-1}}}
%EndExpansion
\dfrac1x\left(  -x+\dfrac{2\gamma^{2j+1}}{\gamma-1}\right)  dx=g(j,\gamma,X)
\end{array}
\]
where
\[
g(j,\gamma,X)=\left(  X+\gamma^{2j-1}\right)  -\frac{2\gamma^{2j+1}}{\gamma
-1}\left(  \ln(-X)-(2j-1)\ln(\gamma)\right)  .
\]
Clearly
\[
I(X)=\frac1{2X}\left(
%TCIMACRO{\dsum \limits_{p=-\infty}^{i-1}}%
%BeginExpansion
{\displaystyle\sum\limits_{p=-\infty}^{i-1}}
%EndExpansion
f(p,\gamma,\gamma^{2p+2})+f(i,\gamma,X)-%
%TCIMACRO{\dsum \limits_{q=-\infty}^{j-1}}%
%BeginExpansion
{\displaystyle\sum\limits_{q=-\infty}^{j-1}}
%EndExpansion
g(q,\gamma,-\gamma^{2q+1})-g(j,\gamma,-X)\right)
\]
and, if $\left\lceil \ln(X)/\ln(\gamma)-1\right\rceil $ is even, then $i=j $.
Upon summation, it can be proved that the minimum and maximum values of $I(X)$
are, respectively,
\[%
\begin{array}
[c]{ccc}%
1+\dfrac{\gamma(\gamma+1)}{(\gamma-1)^{2}}\ln(\gamma), &  & 1+\dfrac
1e\dfrac{\gamma+1}{\gamma-1}\gamma^{\gamma/(\gamma-1)}.
\end{array}
\]
The former quantity is least when $\gamma=5.7041372673...$; the latter
quantity is least when $\gamma=3.2232549401...$. The corresponding mean ratios
are $4.0089813375...$ and $4.8131558458...$. These values together constitute
our solution to the min-mean problem.

An alternative approach is due to Beck \&\ Newman \cite{BN, Gal, AG, KRT,
Alex}. It uses a single random variable $H$, assumed to be uniformly
distributed on the interval $[0,2)$, to sample different logarithmic coils
with rate of increase $\gamma$. For simplicity, take $X>0$. Then
\begin{align*}
\operatorname*{E}(\delta(X))  &  =X+2\operatorname*{E}\left(
%TCIMACRO{\dsum \limits_{j=-\infty}^{2i+1}}%
%BeginExpansion
{\displaystyle\sum\limits_{j=-\infty}^{2i+1}}
%EndExpansion
\gamma^{j+H}\;|\;\gamma^{2i+H}<X\leq\gamma^{2i+2+H}\right) \\
\  &  =X+\frac{2\gamma^{2}}{\gamma-1}\operatorname*{E}\left(  \gamma
^{2i+H}\;|\;\frac X{\gamma^{2}}\leq\gamma^{2i+H}<X\right) \\
\  &  =X+\frac{2\gamma^{2}}{\gamma-1}\operatorname*{E}\left(  X\,\gamma
^{H-2}\right) \\
\  &  =X+\frac{2\gamma^{2}}{\gamma-1}%
%TCIMACRO{\dint \limits_{0}^{2}}%
%BeginExpansion
{\displaystyle\int\limits_{0}^{2}}
%EndExpansion
X\,\gamma^{h-2}\frac12dh=X\left(  1+\frac{\gamma+1}{\ln(\gamma)}\right)
\end{align*}
and this is least when $\gamma=3.591121476669...=1/W(1/e)$, where $W$ denotes
Lambert's function. A search strategy as such is called a mixed strategy (in
game theory)\ or a random strategy (in computer science). Note, however, that
a uniform distribution on $H$ does not imply a uniform distribution on $X$. It
is not clear to us whether the optimal mixed strategy (with $\gamma=1/W(1/e)$)
is necessarily preferable to our deterministic strategy discussed earlier.

\subsection{Acknowledgements}

We grateful to Ricardo Baeza-Yates, John Shonder and Scott Burlington for
their assistance. A discussion during the first author's Visiting Lecture at
Oberlin College in November 2004 was also very helpful.

\end{document}